\documentclass[letter,11pt,round]{article}

\RequirePackage[OT1]{fontenc}
\usepackage{amsthm,amsmath,natbib}
\RequirePackage{hypernat}
\usepackage[pdftex]{graphicx}
\usepackage{amsfonts}
\usepackage{epic}
\usepackage{fancybox}
\usepackage[pdftex]{color}

\usepackage[letterpaper]{geometry}
\usepackage{ulem}

\RequirePackage{ifthen}

\newtheorem{theorem}{Theorem}
\numberwithin{equation}{section}
\theoremstyle{plain}

\newtheorem{remark}{Remark}[section]
\newtheorem{definition}{Definition}[section]

\newtheorem{example}{Example}[section]

\newcommand{\Prob}[1]{\mbox{Pr}\left(#1\right)}

\newcommand{\Pl}[1]{\mbox{Pl}\left(#1\right)}

\newcommand{\Cauchy}[1]{\mbox{Cauchy}\left(#1\right)}
\newcommand{\Bel}[1]{\mbox{Bel}\left(#1\right)}

\definecolor{gray}{rgb}{0.3,0.3,0.3}
\begin{document}
\bibliographystyle{chicago}

%\red

%\begin{frontmatter}
%\title{Point Estimation of Many Normal Means}
%\title{A New Approach to Estimation of the Mean of Multivariate Normal Distribution}
%\title{Many Normal Means: James-Stein Estimators and Inferential Models}
%\title{Elucidating Foundations of Statistical Inference with the Cauchy Distribution}
\title{Elucidating Inferential Models with the Cauchy Distribution}
%\title{Elucidating Inferential Models with the Cauchy Distribution}
%\runtitle{Many Normal Means}

%\begin{aug}
\author{
Chuanhai Liu\footnote{Department of Statistics, Purdue University}\qquad
and\qquad
Ryan Martin\footnote{
Department of Statistics, North Carolina State University}
%  \fnms{Chuanhai} \snm{Liu}
%\thanksref{t2}\ead[label=e1]{chuanhai@purdue.edu}
}

\maketitle

%\thankstext{t2}{Department of Statistics, Purdue University,
%    150 N. University Street, West Lafayette, IN 47907.}%  \printead{e1,u1}.}
%\runauthor{C. Liu}

%\affiliation{Academia Sinica and Purdue University}

%\address{Department of Statistics, Purdue University, 250 N. University Street,
%West Lafayette, IN 47907.  \printead{e1,u1}.}

%\end{aug}

\begin{abstract}
Statistical inference as a formal scientific method  to covert experience to
knowledge  has proven to be elusively difficult.
While frequentist and Bayesian methodologies have been accepted in
the contemporary era as two dominant schools of thought,
it has been a good part of the last hundred years to see
growing interests in development of more sound methods,
both philosophically, in terms of scientific meaning of inference,  
and mathematically, in terms of exactness and efficiency.
These include Fisher's fiducial argument, the Dempster-Shafe theory
of belief functions, generalized fiducial, Confidence Distributions,
and the most recently proposed inferential framework, called
Inferential Models.
Since it is notoriously challenging to make exact and efficient 
inference about the Cauchy distribution, 
this article takes it as an example to elucidate different
schools of thought on statistical inference.
It is shown that the standard approach of Inferential Models
produces exact and efficient prior-free probabilistic inference on the
location and scale parameters of the Cauchy distribution,
whereas all other existing methods suffer from various difficulties.
\end{abstract}

\section{Introduction}

Although statistical inference has been widely used in scientific research,
it remains an elusive problem
in the development of statistical foundations.
A good introduction is the 250-year debate between two contemporary
dominant schools of thought, namely, Bayesian and frequentist. 
Here, we provide a brief relevant discussion
on the two dominant schools
from the point of view of converting experience to knowledge.
In the same vein,
brief comments are also given on 
Fisher's fiducial argument, the Dempster-Shafe theory
of belief functions, generalized fiducial, and Confidence Distributions.

From the point of view of converting experience to knowledge,
\cite{martin2015inferential}
articulates that a sound inferential framework should
produce {\it valid} %{\it meaningful}
and {\it efficient}
uncertainty assessments for assertions of interest
on unknown quantities from known information.
A satisfactory special case is Bayesian inference with 
known prior distributions on everything that is unknown,
since the inferential output is {\it valid}, {\it efficient}, and
{\it probabilistic}.
Nevertheless, this does not necessarily mean that
scientists who make inference using the posterior
distributions in this case are Bayesian.
Bayesian statisticians often refer to people who make Bayesian inference
by specifying prior distributions
for everything unknown either based on their subjective opinions, also called
personal probabilities, or aiming to produce frequency-calibrated
inferential results.  To distinguish these two kinds of Bayesian, 
the latter can perhaps be referred to as
{\it generalized} Bayesian, a term used in the study of point inference based on Bayesian decision theory.
As a result, from the point of view of \cite{martin2015inferential},
the most questionable of such Bayesian methodology is the {\it validity} 
or meaningfulness of its inferential results,
a key point in the centuries-long Bayesian-frequentist debate.

It is arguable that the demand for valid uncertainty assessments in
scientific research with somewhat sacrifice of the situation-specific
nature of desirable probabilistic inference has made frequentist
methodology popular in the twentieth century.
The lack of desirable situation-specific probabilistic interpretation
of confidence intervals and p-values,
the two most commonly used frequentist concepts in practice, is unsettling,
at least philosophically.
The recently developed framework of inferential models (IMs) 
helps solve this problem to some extent. More precisely, 
%Martin and Liu (2013, JASA; 2014, {\it Sinica}) 
\cite{martin2014note}
have shown that
under mild conditions, confidence intervals are plausibility
intervals and p-values are plausibilities,
where the concept of plausibility is situation-specific and probabilistic.
In this sense, IM includes a refined school of frequentist as its
special case.

%{\color{blue}
The Cauchy probability distribution plays an important role
in physics, mathematics, and several related disciplines, such as spectroscopy.
It is also one of two practically extreme cases of
the family of student-t distributions, with the familiar Gaussian distribution
as the other case. Since it is notoriously difficult to make exact 
inference about the Cauchy distribution,
%A brief review of methods on parameter estimation can be
%found in Zhang (2010). % and the Wikipedia article 
%({\tt https://en.wikipedia.org/wiki/Cauchy_distribution}).
we take it as an example to elucidate different
schools of thought on statistical inference.
%}
It is shown that IM methods provide exact and efficient inference about both the location and scale parameters of the Cauchy distribution from a sample, whereas other schools of thought suffer various difficulties.

The Cauchy distribution and its estimation are reviewed in Section \ref{sec:cauchy}. We give a brief review of IM in Section \ref{sec:IM}. The IM procedures are illustrated with its solution to estimate the Cauchy location parameter with known scale from a single observation. Inference about the location parameter with known scale and simultaneous inference about both the location and scale parameters are discussed in Section \ref{sec:CIM}. Marginal inference about the location parameter alone 
and the scale parameter alone is presented in Section \ref{sec:MIM}. The inferential problem in this case provides a benchmark example for evaluating differential schools of thought. Bayesian and fiducial methods for this problem are thus discussed. 
We conclude the paper with a few remarks in Section \ref{sec:conclusion}. It is argued there that although our discussion is limited to inference about the Cauchy distribution, the specific IM methods can be extended to making inference about location and scale parameters in general, that is, group transformations.

\section{The Cauchy Distribution}
\label{sec:cauchy}

The standard Cauchy distribution $C(0,1)$ is given by the
probability density function
\[
     f(x; 0, 1) = \frac{1}{\pi}\frac{1}{1+x^2}, \qquad(x \in \mathbb R)
\]
or, equivalently, by the cumulative distribution function
\[
     F(x; 0, 1) = \frac{1}{\pi}\arctan(x) + \frac{1}{2}, \qquad(x \in \mathbb R)
\]
The general Cauchy distribution $C(\mu, \sigma)$ is obtained
easily from $C(0,1)$ by introducing a location parameter $\mu
\in {\mathbb R}$ and scale parameter $\sigma \in {\mathbb R}_+ = \{\sigma: \sigma \in {\mathbb R}, \sigma>0\}$.
It follows that $C(\mu, \sigma)$ has the probability density function
%and cumulative distribution function
given by
\begin{equation}\label{eq:cauchy-pdf}
f(x; \mu, \sigma) = \frac{1}{\pi\sigma}\frac{1}{1+\frac{(x-\mu)^2}{\sigma^2}},
%\quad \mbox{ and } \quad
%F(x; \mu, \sigma) = \frac{1}{\pi}\arctan(\frac{x-\mu}{\sigma}) + \frac{1}{2},
\qquad(x \in \mathbb R)
\end{equation}
where $\theta \equiv (\mu, \sigma)' \in \Theta \equiv 
{\mathbb R}\otimes {\mathbb R}_+.$
The problem we consider in this paper
is to make inference about the location parameter
$\mu$ with known or unknown scale $\sigma$ from a sample of size $n$,
$X_1, ..., X_n$. 
%, which are {\it iid} $C(\mu, \sigma)$.

As %Zhang (2010) 
\cite{zhang2010highly}
pointed out, it appears notoriously difficult to 
estimate the parameters of $C(\mu, \sigma)$.
%The maximum likelihood estimate of the two-parameter model,
%$\hat\mu$ and $\hat\sigma$, is given by the system of two equations: 
%\begin{eqnarray} \label{eq:MLE-eq01}
%\sum_{i=1}^n \frac{x_i - \hat\theta}{\hat\sigma^2 +(x_i-\hat\mu)^2} &=& 0,\\
%\label{eq:MLE-eq02} \frac{n}{\hat\sigma} - 
%\sum_{i=1}^n \frac{2\hat\sigma}{\hat\sigma^2 +(x_i-\hat\mu)^2} 
%&=& 0.
%\end{eqnarray}
For $n>3$, %Copas (1975) and Gabrielsen (1982) 
\cite{copas1975unimodality} and \cite{gabrielsen1982unimodality}
showed that the joint likelihood function of $\mu$ and $\sigma$ is unimodal.
Although this appears attractive as far as computational methods,
such as Newton-Raphson and EM algorithms, are concerned,
the fact that the likelihood function of $\mu$ with fixed $\sigma$
is randomly multimodal,
with the number of local maxima (excluding the global maxima) asymptotically Poisson distributed with mean $1/\pi$,
is troublesome; %See Barnett (1966), Edwards (1972), and Reeds (1985).
see %Barnett (1966), Reeds (1985), and Bai and Fu (1987).
\cite{barnett1966evaluation}, \cite{reeds1985asymptotic}, and \cite{bai1987maximum}.

The Cauchy distribution is known to be stable, that is,
$\frac{1}{n}\sum_{i=1}^n X_i \sim C(\mu, \sigma)$. This implies that the sample
mean as an estimator of $\mu$ is inefficient.
Alternative methods such as the trimmed mean method and the more general
L-estimators have been investigated in the literature; see
%Rothenberg {\it et al.} (1964), Bloch (1966), Cane (1974),
%and Zhang (2010). 
\cite{rothenberg1964note}, \cite{bloch1966note}, \cite{cane1974linear},
and \cite{zhang2010highly}.
%{\color{red} and {\tt https://en.wikipedia.org/wiki/Cauchy\_distribution}.  }

Another interesting estimator of $\mu$ with fixed $\sigma$ is known as the Pitman estimator %(Pitman, 1939; Cohen Freue, 2007). 
\citep{pitman1939tests,freue2007pitman}.
This estimator is the mean of the Bayesian posterior of $\mu$ computed
with the flat prior, that is,
\begin{equation}\label{eq:pitman-estimator}
\hat{\mu}_P = 
\frac{
\int_{-\infty}^{\infty} u \prod_{i=1}^n f(x_i; u, \sigma) du
}{
\int_{-\infty}^{\infty} \prod_{i=1}^n f(x_i; u, \sigma) du
}
\end{equation}
provided $n>1$.
We show in Section \ref{sec:CIM} that inference using the Bayesian posterior
of $\mu$ with the flat prior is in a certain sense equivalent to  
that produced by the IM method.
This result can be established by applying \cite{lindley1958fiducial} %Lindley (1958)
to a conditional likelihood function. 
The frequentist conditional inference approach is related to
the approach we call conditional IMs.

%McCullagh (1992) 
\cite{mccullagh1992conditional} investigated conditional inference on
Cauchy models with both location and scale parameters.
To make inference on the location parameter $\mu$,
%McCullagh (1992)
\cite{mccullagh1992conditional} took Fisher's fiducial approach
and demonstrated that Fisher's fiducial approach
generates different results for two different ancillary statistics.
This discovery adds to the interesting examples on
the difficulties of Fisher's fiducial method.
More detailed discussion on this problem is given in Section \ref{sec:fiducial-difficulty}.

\section{Inferential Models}
\label{sec:IM}

\subsection{The basic framework}

To produce prior-free probabilistic inference,
the key idea of IM is to produce uncertainty
assessments by predicting unobserved predictable random quantities
based on observed events.
A random variable is called predictable {\it iff}
its distribution is fully known.
%This idea is so obvious that no argument is necessary.
%Because of this simple but fundamental idea, IM automatically includes valid Bayes as a special case.
Let $X$ denote the observed data in its sampling space ${\mathbb X}$,
and let $\theta$ denote the unknown parameter in the parameter space $\Theta$.
%used to specify the sampling distribution of $X$.
Technically, for the postulated model, we introduce
an auxiliary unobserved random quantity $U \in {\mathbb U}$, which has a known distribution,
and write the model using a system of equations that
associate $X$ with $\theta$ and $U$:
\begin{equation}
\label{eq:association-01}
 a(X, \theta, U) = 0
\end{equation}

IM inference is obtained by employing the so-called valid random sets
(see Section \ref{sec:random-set}) to predict the
unobserved realization of $U$, denoted by $U^\star$ that is associated
with the observed data $X$ through \eqref{eq:association-01}. More precisely, let $\mathcal{S}$ denote the predictive random set.
For any assertion
$\theta \in A \subset \Theta$ of interest, the evidence supporting
$A$ is defined as the probability that all possible values of $\theta$
satisfy (\ref{eq:association-01}) 
\begin{equation}
 \Theta_{X, S}  = \left\{
 \theta: a(X, \theta, U)=0 \mbox{ for some } U \in {\cal S}
 \right\}
\end{equation}
fall into $A$, that is,
\begin{equation}
 \Bel{A|X} = \Prob{\Theta_{X, S} \subseteq A}.
\end{equation}
The function $\Bel{.|X}$ is called the belief function.
For a recent discussion on the concept of belief functions,
see \cite{martin2010dempster} %Martin {\it et al.} (2010) 
and reference therein.

Unlike the usual probability, which is additive, {\it i.e.},
$\Prob{A}+\Prob{A^c} = 1$, the belief function is sub-additive, {\it i.e.,}
\[
\Bel{A|X}+\Bel{A^c|X} \le 1
\]
for all $A \subseteq \Theta$.
The fact that the belief function is sub-additive is 
due to the complication that $\Theta_{X, S}$
can overlap with both $A$ and $A^c$.
This overlap creates a new case against neither of $A$ and $A^c$. A complete
characterization of uncertainty assessment of $A$  is then
achieved by introducing one more function, called the plausibility function
that is defined as
\begin{equation}
 \Pl{A|X} = \Prob{\Theta_{X, S} \cap A \neq \emptyset}
 \qquad (A\subseteq \Theta).
\end{equation}
It follows that
\begin{equation}
 \Pl{A|X} =  1 - \Bel{A^c|X}\qquad (A\subseteq \Theta).
\end{equation}

\begin{example}
\label{example:01}
\rm
Consider inference about $\mbox{Cauchy}(\mu, 1)$ from a 
single observation $X$.
Write the sampling model via the following association
\[
X = \mu + Z\qquad (Z\sim \Cauchy{0, 1}).
\]
Suppose that simultaneous inference on all singleton assertions $\{\mu\}$ 
is of interest.
Take the symmetric random set
\[
{\cal S}=[-|C|, |C|]\qquad (C \sim \Cauchy{0,1}).
\]
It generates the following inferential results:
\begin{enumerate}
\item $\Bel{\{\mu\}|X} = 0$; and
\item $\Pl{\{\mu\}|X} = \Prob{\{X-\mu\} \in {\cal S}}
= \Prob{|X-\mu| \leq |C|}
= 2\Prob{C \leq -|X-\mu|},
%= 1 - \Prob{|C| \leq |X-\mu|},
$
%given by the {\it cdf} of $|C|$,
\end{enumerate}
for all $\mu \in (-\infty, \infty)$.

\begin{center}
%a quantile-quantile plot
\vspace{-0.1in}
\begin{figure}[!htb]
%\begin{center}
%\subfigure{\scalebox{0.30}{\includegraphics{Programs/2D-naive.pdf}}}
%\includegraphics[height=3cm, bb=0 0 100 100]{Programs/2D-naive.pdf}
%\scalebox{0.30}{\includegraphics{Programs/2D-naive.pdf}}
%\includegraphics[height=12cm, bb=0 0 100 100]{plots/plau-func-01.pdf}
%\includegraphics[width=0.95\linewidth]{plots/plau-func-01.pdf}
\includegraphics[width=0.95\linewidth]{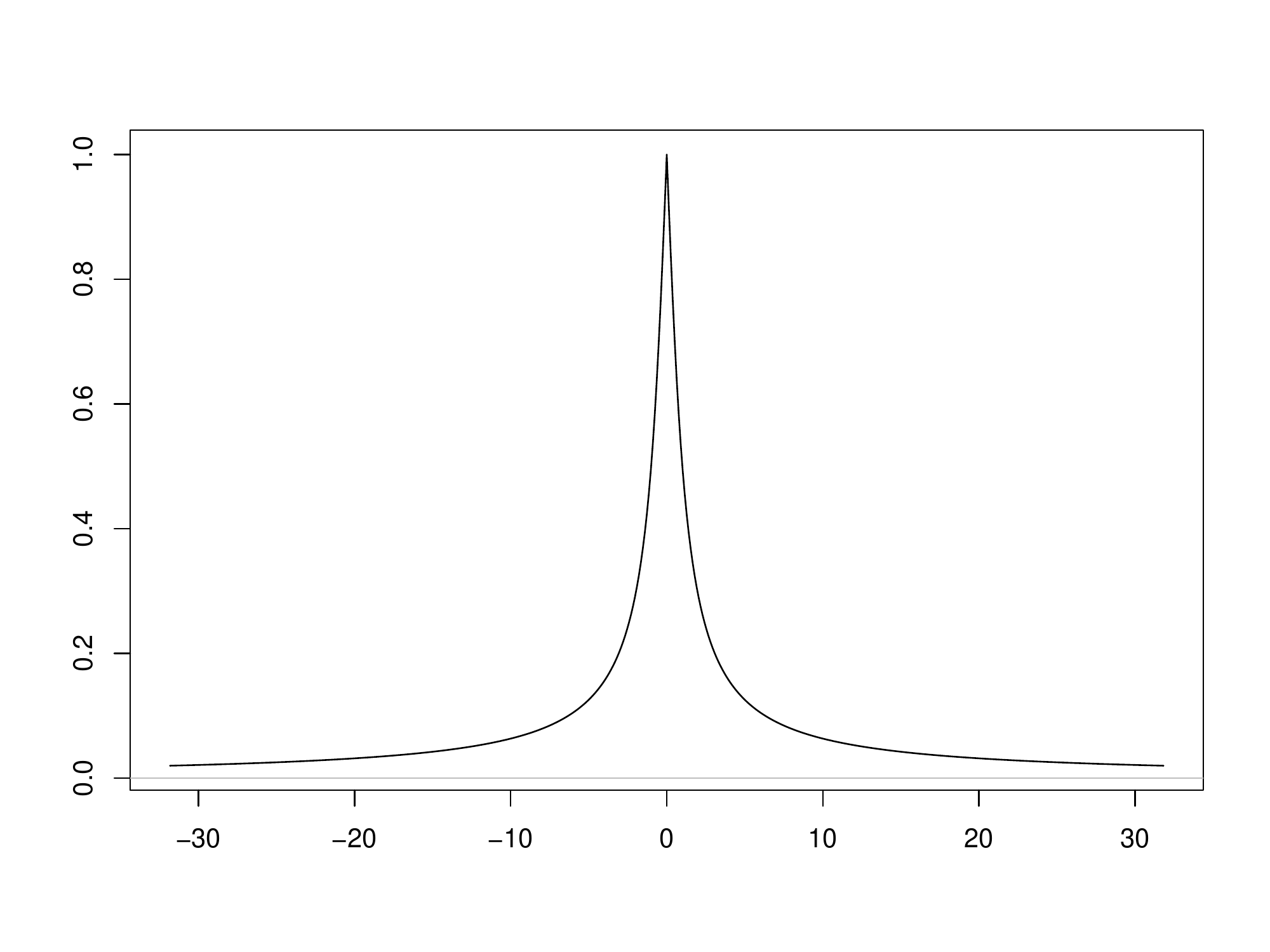}
%\hspace{-1.2in}
%\subfigure{\scalebox{0.31}{\includegraphics{IM_02.pdf}}}
%\caption{The first two pages of a 70-page draft lecture notes by Prof N. Hung.}
\caption{
The plausibility function $\Pl{\{\mu\}|X=0}$
in Example \ref{example:01}.
}
\label{fig:plau-func-01}
%\end{center}
\end{figure}
\end{center}

The plausibility curve given $X=0$ is depicted in Figure
\ref{fig:plau-func-01}.
For a point estimate of $\mu$, one may use the most plausible
value of $\mu$. The most plausible value depends on the 
construction of the underlying random set.
However, it is seen from this example that all the uncertainty
is characterized by the distribution of the predictable random variable
$Z$, which is known as the pivotal quantity in the statistics literature.
The use of plausibility on assertions on $\mu$ instead of assigning the usual probability to $\mu$ is to avoid the difficulty
in inducing the probability for $\mu$ through that of $Z$ via
the association $\mu = X - Z$.
More discussion of this issue is given in Section
\ref{sec:fiducial-difficulty}.

\end{example}

\subsection{Random sets}
\label{sec:random-set}

%\begin{remark}[IM versus fiducial] \rm
Although the idea of IM was motivated by Fisher's fiducial argument
for prior-free probabilistic inference 
\citep{martin2010dempster,zhang2011dempster,martin2013inferential},
IM reasons with uncertainty by working with the auxiliary random variable $U$.
Since $U$ is well defined, nothing can go wrong in probability operations.
This is not the case for Fisher's fiducial argument. Assume, for example, that
there exists a one-to-one mapping between $U$ and $\theta$ for given $X$.
Write
\begin{equation}
\label{eq:fiducial-map-01}
\theta = \theta_X(U).
\end{equation}
Unfortunately, the mapping (\ref{eq:fiducial-map-01})
does not allow for valid probability mass propagation, as the
mapping changes from one pair of $(U, \theta)$ to another.
This explains why IM uses random sets to predict
$U^\star$. The use of a random set effectively carries out
probability calculations/integration.
%\end{remark}
This intuition behind IM helps to establish the desired frequency-calibration property, which is termed validity in the IM framework.

\begin{definition}[Validity]
A random set ${\cal S}\subset {\mathbb U}$, independent of $U^\star$,
is said to be valid
to predict $U^\star$ iff the probability that ${\cal S}$
covers $U^\star$ is not smaller in distribution
than the standard uniform random variable, that is,
\begin{equation}
\Prob{\Prob{U^\star \in {\cal S}|U^\star} \le u} \ge u
\end{equation}
for all $u \in[0, 1]$.
\end{definition}

%Martin and Liu (2013) 
\cite{martin2013inferential} shows that
the use of valid random sets guarantees that the resulting inference about
$\theta$ is valid in the sense that
under the truth of $A$, the belief function $\Bel{A^c|X}$ as a measure of 
evidence against $A$ is bounded from above by the standard uniform random random
variable in repeated experiments.

\section{Combining information}
\label{sec:CIM}

\subsection{Inference on $\mu$ with known $\sigma$}
\label{subsec:CIM-4-mu}

%\begin{example} \rm
Consider inference on the Cauchy distribution
with unknown location $\mu$ and known scale $\sigma$ from a sample size of
$n$, $X_1, ..., X_n$.  Without loss of generality, we take the unit scale $\sigma=1$.
A simple association is given by
\begin{equation}\label{eq:CIM-A-01}
X_i = \mu + U_i\qquad(i=1, ..., n; U_i\stackrel{iid}{\sim} \mbox{C}(0, 1))
\end{equation}
%\end{example}
For valid inference about $\mu$, the idea of predicting
unobserved $U_i$' with predictive random sets in the $n$-dimensional space ${\mathbb R}^n$
still works but is inefficient. For efficient inference via prediction,
we rewrite (\ref{eq:CIM-A-01}) as
\begin{eqnarray}
X_1 &=& \mu + U_1\label{eq:CIM-A-02-01} \\
W_i \equiv U_i - U_{1} &=&X_i - X_{1} \qquad(i=2, ..., n)
\label{eq:CIM-A-02-02}
\end{eqnarray}
It is clear that inference about $\mu$ can be made by predicting 
$U_1$ alone. Furthermore, the information given by
the $n$-dimensional quantity $W\equiv(W_2, ..., W_n)'$ in (\ref{eq:CIM-A-02-02})
can be used to make a more efficient prediction of $U_1$.
This idea of making efficient predictions motivated what is called conditional IM. General CIM theory and its relationship with
Bayesian inference and Fisher conditional inference are established by 
\cite{martin2015conditional}. For example, \cite{martin2015conditional} showed that CIM procedures produce valid inference.

The conditional distribution of $U_1$ can be obtained straightforwardly
with routine algebraic operations and is given by
%\ref{eq:cauchy-pdf}
\begin{equation}
\label{eq:CIM-U1-pdf}
f_{U_1|W}(u|W=w)
\propto
%\frac{1}{\pi}
\frac{1}{1+u^2} \prod_{i=2}^n % \frac{1}{\pi}
\frac{1}{1+(u+w_i)^2}
\qquad(u \in {\mathbb R})
\end{equation}
IM belief and plausibility values can be evaluated accordingly,
using standard numerical methods.
%Martin and Liu (2013) 

It appears that the CIM construction given by
(\ref{eq:CIM-A-02-01}) and (\ref{eq:CIM-A-02-02}) is not unique. 
The key fact is that a $(n-1)$-dimensional quantity of
$U\equiv(U_1,...,U_n)'$ is fully observed.
To make inference about $\mu$, we only need to predict a one-dimensional
quantity given the observed information.
The following result addresses the issue of uniqueness 
in choosing different constructions from the following
class of CIM index by
$a \equiv(a_1,...,a_n)' \in {\mathbb R}^n$ with $\sum_{i=1}^n a_i = 1$
and $a_1 \neq 0:$
\begin{eqnarray}
\sum_{i=1}^n a_iX_i &=& \mu + T\\
\label{eq:CIM-GA-02-01}
W_i \equiv U_i - T &=&X_i - \sum_{j=1}^n a_jX_j \qquad(i=2, ..., n)
\label{eq:CIM-GA-02-02}
\end{eqnarray}
where
$T \equiv \sum_{i=1}^n a_iU_i$.

\begin{theorem}
\label{thm:CIM-post}
For any predictive random set on $T$, ${\cal T}$, conditional on $W=w$, the set 
\begin{equation}
\{\mu: \mbox{$\sum_{i=1}^n a_iX_i = \mu + t$ for some $t\in {\cal T}$} \}
\end{equation}
is the same for all $a$, where $a \equiv(a_1,...,a_n)' \in {\mathbb R}^n$ with $\sum_{i=1}^n a_i = 1$
and $a_1 \neq 0$.
\end{theorem}

\begin{proof}
It is clear from (\ref{eq:CIM-GA-02-01}) that we need only to show
that the conditional probability that $T \geq \sum_{i=1}^n a_iX_i - \mu_0$ given
$W=w$ is independent of $a$ for all $\mu_0 \in {\mathbb R}$.
Due to the fact that the Jacobian of the transformation from $(U_1,..., U_n)'$ to $(T, W)$ is $|a_1|$, we have
\begin{eqnarray*}
f_{T|W=w}(t) &\propto& \frac{1}{1+ 
\frac{(t-\sum_{k=2}^na_k (t+X_k - \sum_{j=1}^n a_j X_j))^2}{a_1^2}}
\prod_{k=2}^n \frac{1}{1+ (t + X_k-\sum_{j=1}^n a_j X_j)^2}\\
&=& 
\prod_{k=1}^n \frac{1}{1+ (t + X_k-\sum_{j=1}^n a_j X_j)^2}\\
\end{eqnarray*}
for $t \in {\mathbb R}$.
This leads to 
\begin{eqnarray}
\Prob{ T \geq \sum_{i=1}^n a_iX_i - \mu_0 } &\propto&
=\frac{
\int_{-\infty}^{\mu_0} \prod_{k=2}^n \frac{1}{1+ (X_k - u)^2} du
}{
\int_{-\infty}^{\infty} \prod_{k=2}^n \frac{1}{1+ (X_k - u)^2} du
}
\label{eq:CIM-posterior}
\end{eqnarray}
for all $t_0 \in {\mathbb R}$, completing the proof.

It should be noted that although a more general class of CIM
that represents the $(n-1)$-dimensional quantity with
an arbitrary linear transformation of $(U_1-T,..., U_n-T)$
can be considered, it is easy to see that the conclusion
remains the same, except that the corresponding argument needs
to use more tedious notation.
\end{proof}

\begin{remark}
It is seen from the result (\ref{eq:CIM-posterior})
in the proof of Theorem \ref{thm:CIM-post} that
efficient inference about the location parameter $\mu$ with known $\sigma$
is intrinsically the same as that offered by Bayesian posterior
with flat prior on $\mu$.
This conclusion holds for any sample $X_1, ..., X_n$ from the joint
distribution of the form $F(X_1 -\mu, ..., X_n-\mu)$.
This result can be established by applying Lindley (1958)
to the likelihood of $\mu$ given $X_1$ but from the conditional
model of $X_1$ given $W$.
\end{remark}

\subsection{Joint inference on $\mu$ and $\sigma$}

For joint inference about $\mu$ and $\sigma$, the association for the sample of size $n$ $(\ge 2)$, $X_1,..., X_n$,
from $\mbox{C}(\mu, \sigma)$ can be written as
\begin{equation}\label{eq:CIM-A-04}
X_i = \mu + \sigma U_i\qquad(i=1, ..., n; U_i\stackrel{iid}{\sim} \mbox{C}(0, 1))
\end{equation}
Following \cite{fisher1934two} and \cite{mccullagh1992conditional}, %McCullagh (1992), 
we
use the order statistic 
$(X_{(1)}, ..., X_{(n)})'$ and write an association for combining information as
\begin{eqnarray}
% \frac{X_{(1)}-\mu}{X_{(2)}-X_{(1)}} &=& \frac{T}{S}
 X_{(1)} &=& \mu + \sigma T
 \label{eq:CIM-A-04-01}
 \\
 X_{(2)} - X_{(1)} &=& \sigma S
 \label{eq:CIM-A-04-02}
 \\
 \frac{X_{(i)} - X_{(1)}}{X_{(2)}-X_{(1)}} &=&  W_{i}, \qquad(i=3,...,n)
 \label{eq:CIM-A-04-03}
\end{eqnarray}
where $(U_{(1)}, ..., U_{(n)})'$ denotes that ordered
$U_1, ..., U_n$  with $U_{(1)} \leq ...\leq U_{(n)}$,
$T = U_{(1)},$ $S = U_{(2)}- U_{(1)},$ and
$ W_i =  \frac{U_{(i)} - U_{(1)}}{U_{(2)}-U_{(1)}} $ for $i=3,...,n.$
Thus, inference about $(\mu, \sigma)$ can be obtained by
predicting $(T, S)$ from its predictive distribution given
the observed information $(W_3, ..., W_n)$.

\begin{theorem}
The conditional density of $(T, S)$ given
$W=(W_3,...,W_n)' = (w_3,..., w_n)$ is given by 
\begin{equation}\label{eq:CIM-TS-pdf}
f_{T, S | W=w}(t, s) \propto
s^{n-2}
\frac{1}{1+t^2}
\frac{1}{1+(t+s)^2}
\prod_{i=3}^n \frac{1}{1+(t+w_is)^2},
\qquad(t\in {\mathbb R}, s \in {\mathbb R}_+)
\end{equation}
\end{theorem}

\begin{proof}
It is known that the density function of $(U_{(1)}, ..., U_{(n)})'$ 
is given by
\[
f_{U_{(1)},...,U_{(n)}}(u_1, ..., u_n) =
n! \prod_{i=1}^n \frac{1}{\pi}\frac{1}{1+u_i^2},
  \qquad (u_1 \leq u_2 \leq ... \leq u_n).
\]
The Jacobean of the transformation
\begin{eqnarray*}
U_{(1)} &=& T\\
U_{(2)} &=& T + S\\
U_{(i)} & = & T + S W_i, \qquad(i=3,...,n)
\end{eqnarray*}
can be easily obtained as $s^{n-2}$.
It follows that the conditional density of $(T, S)$ given
$W=(W_3,...,W_n)' = (w_3,..., w_n)$ is given by 
\[
f_{T, S | W=w}(t, s) \propto s^{n-2} \frac{1}{1+t^2} \frac{1}{1+(t+s)^2}
\prod_{i=3}^n \frac{1}{1+(t+w_is)^2},
\qquad(t\in {\mathbb R}, s \in {\mathbb R}_+)
\]
\end{proof}

%\label{subsec:CIM-4-mu}

\section{Marginal inference}
\label{sec:MIM}

In this section, we discuss inference on $\sigma$ alone
and inference on $\mu$ alone in the two-parameter Cauchy model $C(\mu, \sigma)$.
Both of the problems belong to
so-called marginal inference, inference on lower-dimensional quantity of
unknown parameters. %Martin and Liu (2015) 
\cite{martin2015marginal} discussed a general approach to
marginal inference in the IM framework. 
It turns out that the IM inference on $\sigma$ with unknown $\mu$ and
inference on $\sigma$ with unknown $\mu$ are consistent with Bayesian
inference with the Pitman prior
\begin{equation}
\label{eq:pitman-prior}
  \frac{1}{\sigma} d\mu d\sigma, \qquad(\mu \in {\mathbb R}, \sigma \in {\mathbb R}_+)
\end{equation}
in terms of Bayesian credible intervals and IM plausibility intervals.

%McCullagh (1992) 
\cite{mccullagh1992conditional} showed that the two-parameter Cauchy model provides
an example of non-uniqueness of Fisher's fiducial inference, based on
two different ancillary statistics. Section \ref{subsec:mim-03} discusses
the relevant
issues from the perspective of marginal IM.

\subsection{Inference on $\sigma$ with unknown $\mu$}
\label{subsec:mim-01}

The basic idea to make efficient marginal inference on $\sigma$ is to
produce valid inference by predicting a one-dimensional auxiliary random
variable rather than two-dimensional auxiliary random variables.
An intuitive explanation of its efficiency is as follows.
For any valid predictive random set ${\cal S}_{T,S}$ for the unobserved $(T, S)$
defined in (\ref{eq:CIM-A-04-01}) and (\ref{eq:CIM-A-04-02}), the projected
predicted random set 
\[
{\cal S}_S = \{s: (t, s) \in {\cal S}_{T,S} \mbox{ for some } t\}
\]
provides a valid predictive random set for $S$ alone and, thereby,
for valid marginal inference on $\sigma$ via (\ref{eq:CIM-A-04-01}).
This implies that efficient inference on $\sigma$ can be obtained
by creating valid two-dimensional predictive random sets 
that have stochastically small projections on the space of $S$. 

For the two-parameter Cauchy mode,
efficient marginal inference about $\sigma$ can be obtained
by predicting the unobserved $S$ in (\ref{eq:CIM-A-04-02}) from
(\ref{eq:CIM-TS-pdf}) or its
marginal distribution
\[
f_{S| W=w}(s) \propto s^{n-2} \int_{-\infty}^\infty \frac{1}{1+t^2} \frac{1}{1+(t+s)^2}
\prod_{i=3}^n \frac{1}{1+(t+w_is)^2} dt,
\qquad(s \in {\mathbb R}_+)
\]
Let ${\cal S}$ be a valid predictive random set for the unobserved 
$S$. The propagated random set of ${\cal S}$  on the space of $\sigma$
is given by
\[
\{(X_{(2)}-X_{(1)})/s: s \in {\cal S}\}.
\]

For the one-sided valid predictive random set
\[
{\cal S}_S = \{s: s \geq S\}, \qquad(S\sim f_{S|W=w}(s))
\]
the propagated random set on the space of $\sigma$ is 
characterized by its upper unbound
\[
     G \equiv \frac{X_{(2)}- X_{(1)}}{S}. 
\]
It is straightforward to show that the density of the $G$ 
has the form of
\begin{equation}
\label{eq:pitman-posterior-sigma}
f_{G| W=w}(g) \propto  \frac{1}{g^{n+1}} \int_{-\infty}^\infty
\prod_{i=1}^n \frac{1}{1+\frac{(X_{(i)}-u)^2}{g^2}} du,
\qquad(s \in {\mathbb R}_+)
\end{equation}
This establishes a connection between IM and the Bayesian approach to
marginal inference on $\sigma$ in the two-parameter Cauchy model
$C(\mu, \sigma)$.

\subsection{Inference on $\mu$ with unknown $\sigma$}
\label{subsec:mim-02}
Similar to marginal inference on $\sigma$,
efficient inference about $\mu$ is obtained by
making use of the marginal association
\begin{equation}
\label{eq:MIM-A-mu}
\frac{X_{(1)}-\mu}{X_{(2)} - X_{(1)}} = \frac{T}{S}
\end{equation}
via prediction of $M \equiv T/S$ conditional on observed $W$.
The density of the predictive distribution for $M$ is given by
\[
f_{M| W=w}(m) \propto  s^{n-2+1} \int_{0}^\infty
\frac{1}{1+m^2s^2} \frac{1}{1+(ms+s)^2}
\prod_{i=3}^n \frac{1}{1+(ms+\frac{X_{(i)}-X_{(1)}}{X_{(2)}-X_{(1)}}s)^2} ds,
\qquad(m \in {\mathbb R})
\]

The above predictive distribution is intrinsically related to Bayesian
inference with the Pitman prior (\ref{eq:pitman-prior}), rather than
the Jeffreys prior. To see this, note that for any predictive random set $S$ on $M$, its propagated random set in the $\mu$ space through
(\ref{eq:MIM-A-mu}), the mapping
\[ \mu = X_{(1)} - (X_{(2)}-X_{(1)}) M
\]
For the case with a one-sided predictive random set for 
$M$, the random set on $\mu$ is also one-sided. The density
of its non-trivial end point is given by
\[
f_{Z| W=w}(z) \propto \int_{0}^\infty \frac{1}{s^{n}}
\prod_{i=1}^n \frac{1}{1+\frac{(z-X_{(i)})^2}{s^2}} \frac{1}{s} ds,
\qquad(z \in {\mathbb R})
\]
This result corresponds to the claimed Bayesian interpretation of the 
IM marginal inference on $\mu$ in the two-parameter Cauchy model.

\subsection{IM versus Fisher's fiducial}
\label{subsec:mim-03}
\label{sec:fiducial-difficulty}

The work on the Inference Models was motivated originally
by fiducial inference and the Dempster-Shafe theory of belief functions.
As articulated in %Liu and Martin (2014), 
\cite{liu2015frameworks}, while IM provides exact prior-free
inference, fiducial inference is not prior-free.
More discussions can be found in 
%Liu and Martin (2014) and Martin and Liu (2015c).
\cite{liu2015frameworks} and \cite{martin2015inferential}.
Here we take the non-uniqueness problem discovered by
%McCullagh (1992) 
\cite{mccullagh1992conditional} on conditional and fiducial inference 
with the two-parameter model as an example to illustrate
the difficulty of the fiducial argument
from our view point of the IM framework.

Based on nice properties of the Cauchy model $C(\mu, \sigma)$ in
the literature, including 
%Menon (1962, 1966),  Pitman and Williams (1967), and Williams (1969),
\cite{menon1962characterization,menon1966another},
\cite{pitman1967cauchy}, and \cite{williams1969cauchy},
%McCullagh (1992) 
\cite{mccullagh1992conditional} introduced an attractive representation of
the parameter $(\mu, \sigma)$ using the complex number 
\[
\theta = \mu + i \sigma
\]
and write $C(\mu, \sigma)$ as $C(\theta)$.
The result of %McCullagh (1992)
\cite{mccullagh1992conditional} that is to be used here can be
summarized into the following theorem.

\begin{theorem}\label{thm:transformation}
If $Y \sim C(\theta)$, then 
\[
 Y^* = \frac{a Y + b}{cY + d} \sim C\left( \frac{a \theta + b}{c\theta + d} \right).
\]
\end{theorem}

Taking $c=1$ and $a=b=d=0$, we see that
\[
   R_i = \frac{1}{X_i}, \qquad (i=1,...,n)
\]
form a sample of size $n$ from $C(1/\theta)$ or
$C\left(
\frac{\mu}{\mu^2+\sigma^2},
\frac{\sigma}{\mu^2+\sigma^2}
\right)$.
%\footnote{ \color{cyan} \[ \frac{1}{\mu+i\sigma} = \frac{\mu-i\sigma}{\mu^2+\sigma^2} = \]}.
The discussion in Sections \ref{sec:CIM} and \ref{sec:MIM}
suggests that efficient exact inference about 
\[
\mu^* \equiv \frac{\mu}{\mu^2+\sigma^2}
\qquad \mbox{ and } \qquad
\sigma^*\equiv \frac{\sigma}{\mu^2+\sigma^2}
\]
be made with the conditional association
\begin{eqnarray*}
R_{(1)} &=& \mu^* + \sigma^*U_{(1)}\\
R_{(2)} - R_{(1)} &=& \sigma^* (U_{(2)}-U_{(1)})\\
\frac{R_{(i)} - R_{(1)}}{R_{(2)} - R_{(1)}} &=&
\frac{U_{(i)}-U_{(1)}}{U_{(2)}-U_{(1)}}, \qquad (i=3,...,n)\\
\end{eqnarray*}
In this case, the corresponding Bayesian interpretation
of this result is to use the Pitman prior for
$(\mu^*, \sigma^*)$:
\[
\frac{1}{\sigma^*} d\mu^* d\sigma^*
= \frac{1}{\sigma^*} d\mu^* d\sigma^*,
\qquad (\mu^* \in {\mathbb R}, \sigma^* \in {\mathbb R}_+)
\]
Fiducial is ambitious in the sense that it takes
this prior effectively as the implicit prior distribution
on $(\mu^*, \sigma^*)$ and leaves an impression that
the resulting inference is prior-free.
As a result, the use of the two different prior distributions
in the two constructions for conditional and fiducial inference
on $C(\mu, \sigma)$,
namely,
\[
\frac{1}{\sigma} d\mu d\sigma
\]
and
\[
\frac{1}{\sigma^*} d\mu^* d\sigma^*
= \frac{\mu^2+\sigma^2}{\sigma} d\mu d\sigma,
\]
leads to conflicting probabilistic inference.

IM doesn't suffer from the difficulty that the fiducial argument has. The transformation in Theorem \ref{thm:transformation} is used only in the IM frame when the problem of interest is inference on $\mu^* = \mu/(\sum^+\sigma^2)$ or on
$\sigma^* = \sigma/(\mu^2+\sigma^2)$. 
Nevertheless, this example provides as a good example for the argument on the
difficulty made in \cite{liu2015frameworks}.
% More elaboration ...

\section{Concluding Remarks}
\label{sec:conclusion}

Statistical inference is no doubt important in scientific investigations.
Development of solid foundations for inference  
has proven to be one of the most difficult problems.
This article elucidates the IM framework with the Cauchy model.
Due to the limited space, the discussion has be been
focused on intuition-based explanations of Basic IM, Conditional IM,
and Marginal IM with the Cauchy model, for which
inference appears to be difficult in all previous schools of
thought. This is evidential by the impressive list of work on the Cauchy model,
including 
\cite{pitman1939tests},
\cite{menon1962characterization,menon1966another},
\cite{rothenberg1964note},
\cite{bloch1966note},
\cite{pitman1967cauchy},
\cite{williams1969cauchy},
\cite{cane1974linear},
\cite{copas1975unimodality},
\cite{gabrielsen1982unimodality},
\cite{reeds1985asymptotic},
\cite{bai1987maximum}, 
\cite{mccullagh1992conditional},
\cite{freue2007pitman},
\cite{zhang2010highly},
\cite{martin2015conditional},
and the references therein.
The popular method of moments and method of maximum likelihood
are standard textbook frequentist approaches. These methods
are inexact and inefficient.
The non-uniqueness problem discovered by %McCullagh (1992) 
\cite{mccullagh1992conditional} shows the difficulty
of the fiducial argument. It also serves as an example that to suggest
care must be taken with Bayesian inference, especially if 
frequency-calibrated inference is desirable.
 
 Although our discussion is limited to inference about the Cauchy distribution, the specific IM methods can be extended to making inference about location and scale parameters in general, that is, group transformations. The relevant topic deserves a separate in-depth investigation.
 
The current progress on IM can be found in %Martin and Liu (2015c)
\cite{martin2015inferential} with new developments, including 
\cite{liu2020inferential} and
\cite{cella2022validity}, to name a few.
While further development of IM is necessary, we believe that
the underlying idea of predicting unobserved random variables
conditional on observed will continue to be fundamental for
development of exact and efficient probabilistic inference.

%\section*{Acknowledgment}

\bibliographystyle{chicago}
\bibliography{mnms}

\end{document}